\documentclass{amsart}
\usepackage{amsmath}

\newcommand{\al}{\alpha}
\newcommand{\be}{\beta}
\newcommand{\nn}{\mathbb{N}}      
 
\newtheorem{theorem}{Theorem}[section]   
\newtheorem{cor}[theorem]{Corollary}     
\newtheorem{lemma}[theorem]{Lemma}         

\theoremstyle{definition}
\newtheorem{definition}[theorem]{Definition}   

\theoremstyle{remark}
\newtheorem{remark}[theorem]{Remark}        

\input epsf.tex
\begin{document}

\title[Moreno-Soc{\'\i}as Conjecture]  
{Generic ideals and Moreno-Soc{\'\i}as conjecture} 

\author[Aguirre et. al.]{Edith Aguirre  \and Abdul Salam Jarrah \and Reinhard
Laubenbacher \and Juan Ariel Ortiz-Navarro \and Rogelio Torrez}
\address [Edith Aguirre]{Department of Mathematics\\ 
University of Texas at El Paso \\
El Paso, TX 79907}
\email{emaguirre@miners.utep.edu}

\address[Abdul Salam Jarrah]{Department of Mathematical Sciences\\
New Mexico State University\\
Las Cruces, NM 88003}    
\email{ajarrah@nmsu.edu} 
 
\thanks{All correspondence should be addressed to the second author.
This work was done during the 2000 Summer Institute in
Mathematics for Undergraduates (SIMU)
at the University of Puerto Rico--Humacao. The authors thank the organizers,
Herbert Medina and Ivelisse Rubio, for their
wonderful work. The second author wishes to thank Keith Pardue
for helpful suggestions and discussions.
Finally, the authors thank the referees for their useful suggestions and
comments.
}

\address[Reinhard Laubenbacher]
{Department of Mathematical Sciences\\
 New Mexico State University\\
Las Cruces, NM 88003}
\email{reinhard@nmsu.edu}

\address[Juan Ariel Ortiz-Navarro]
{Department of Mathematics\\ 
University of Puerto-Rico\\
Humacao, PR 00791}
\email{mr\_ariel@hotmail.com}

\address[ Rogelio Torrez]
{5646 W Heatherbrae Dr.\\ 
Phoenix, AZ 85031}
\email{torrez@asu.edu}

\date{April 3, 2001}
     
\begin{abstract}
Let $f_1, \dots, f_n$ be homogeneous polynomials generating
a generic ideal $I$ in the ring of polynomials in $n$ variables
over an infinite field. 
Moreno-Soc\'ias conjectured that
for the graded reverse lexicographic term ordering, the initial
ideal ${\rm in}(I)$
is a weakly reverse lexicographic ideal. 
This paper contains a new proof of Moreno-Soc{\'\i}as' conjecture for
the case $n=2$.
\end{abstract}

\maketitle

\section{Introduction}
\noindent We begin by introducing the
definitions necessary to understand the conjecture
under study.

\smallskip 
\noindent Let $R=K[x_1, \dots,x_n]=K[{\bf x}]$ be the polynomial 
ring in $n$ variables over an 
infinite field $K$, which is an extension of a base field $F$.  Then
$R$ is an $\nn$-graded ring such that for each $d \in \nn$, $R_d$ is the 
set of all homogeneous polynomials (forms) of total degree equal to $d$,
where $\deg x_1= \dots =\deg x_n = 1$.
 
\smallskip 
\begin{definition}
Let $f = \sum_m \al_m m$ be a form of degree $d$ such that $m$ runs over all 
monomials of degree $d$ and $\al_m \in K$.  Then $f$ is called {\it generic} if 
the coefficients $\al_m$ satisfy the following two conditions: 
\begin{enumerate}
\item If $m \neq m'$, then  $\al_m \neq \al_{m'}$. 
\item The set of all $\al_m$ is algebraically independent over $F$.
\end{enumerate}  
\end{definition}

\smallskip 
\begin{definition}
An ideal $I \subset R$ is called {\it generic} if $I$
can be generated by generic forms $f_1, \ldots, f_r$, 
where \\
$f_i = \sum_{m_{f_i}} \al_{m_{f_i}} m_{f_i}$,   
satisfying the following conditions:
\begin{enumerate}
\item $\al_{m_{f_i}} \neq \al_{m_{f_j}}$ for $i \neq j$. 
\item The set $\{\al_{m_{f_i}} \; : \; 1 \leq i \leq r \}$ is algebraically 
independent over $F$.
\end{enumerate}  
\end{definition}

\smallskip 
\begin{definition}[Graded reverse lexicographic order]
 Let 
${\bf x}^\al,\: {\bf x}^\be \in R$.  Then let
${\bf x}^\be > {\bf x}^\al$ if $\deg {\bf x}^\be > \deg {\bf x}^\al$ or
$\deg {\bf x}^\al = \deg {\bf x}^\be$ and the right-most nonzero 
entry of $\be -\al$ is negative.
\end{definition}

\smallskip 
\begin{definition}
A monomial ideal $J \subset R$ is {\it weakly reverse lexicographic} 
if, whenever ${\bf x}^\be$ is a minimal generator of $J$,
then every monomial of the same degree which 
precedes ${\bf x}^\be$ in the reverse lexicographic term order
must be in $J$.  And $J$ is {\it reverse lexicographic}
if, whenever ${\bf x}^\be \in J$,
then every monomial of the same degree which 
precedes ${\bf x}^\be$ must be in $J$.
\end{definition}

\smallskip 
\begin{remark}
Every  {\it reverse lexicographic} ideal is {\it weakly reverse lexicographic}.
\end{remark}

\smallskip 
\begin{definition}
Let $f=\sum_{\al\in A}a_\al x^\al$ be a nonzero polynomial in $R$
such that $a_\al \neq 0$ for all $\al \in A$. Let $<$ be a term order.
The {\it leading monomial} of $f$ is LM$_{<}(f)=x^\beta$ if $\beta\in A$ and
$x^\al < x^\be$ for all $\al \in A,\al \ne \be$.
\end{definition}

\smallskip 
\begin{definition}
For an ideal $I \subset R$, the {\it initial ideal} of $I$ is the ideal
$in(I)=\langle LM(f) \; : \; f \in I \rangle$.
\end{definition} 

\smallskip 
\noindent We are now ready to state the
Moreno-Soc{\'\i}as Conjecture and some known facts about it.

\medskip 
\noindent {\bf Conjecture} (Moreno-Soc{\'\i}as \cite{MS, P}). 
Let $d_1, \dots, d_n \in \nn$
and $I \subset R$ be a generic ideal generated by a sequence of polynomials
$f_1, \dots, f_n$ of degrees $d_1, \dots, d_n \in \nn$ , and $J= in(I)$,
the initial ideal of $I$ with respect to the graded reverse lexicographic order.
Then $J$ is weakly reverse lexicographic.

\smallskip 
\noindent As pointed out in \cite{P}, this conjecture 
implies many other
conjectures, in particular, Fr\"oberg's Conjecture \cite{F}, 
which gives a formula
for the Hilbert series of generic ideals. 

\smallskip
The Moreno-Soc{\'\i}as Conjecture
is trivial when $n = 1$.  It was proven by Moreno-Soc{\'\i}as for $n=2$ in
his thesis \cite{MS}, as was pointed out to us by one of the referees.
We learnt about the conjecture in \cite{P}, and this paper contains a
proof for the case $n=2$, which is quite different from the one in
\cite{MS}.  Our proof is quite elementary, and was discovered through
extensive calculations using the computer algebra systems MAPLE
and {\it Singular}.  Unfortunately, it seems unlikely that our methods
can be extended to deal with the case of more variables.  Already for
$n=3$, computer calculations become impractical.

\medskip
\section{The two variable case}
\noindent In this section we show that Moreno-Soc{\'\i}as Conjecture 
is true for
the case of two generic forms in the ring of polynomials in two
variables.
Let $R=K[x,y]$ be  the polynomial ring in 2 variables over an 
arbitrary infinite field $K$, with a base field $F$. 

\smallskip 
\noindent Let $n \leq m, \mu = m-n$, and let 
\begin{eqnarray}\label{f1}
f_1 &=& a_{1,1}x^n+ a_{1,2}x^{n-1}y+ \cdots+ a_{1,n}xy^{n-1}+ a_{1,n+1}y^{n}, \nonumber \\
f_2 &=& b_{2,1}x^m+ b_{2,2}x^{m-1}y+ \cdots+ b_{2,m}xy^{m-1}+ b_{2,m+1}y^{m}
\end{eqnarray}
be generic forms of degree $n$ and $m$, respectively, generating the 
generic ideal $I=\langle f_1,f_2 \rangle$ in $R$. Since $n \leq m$, 
one can divide $f_2$ by $f_1$:
$$
f_2 = qf_1 + r,
$$
where
\begin{equation*}
r = a_{2,1}x^{n-1}y^{\mu+1}+ a_{2,2}x^{n-2}y^{\mu+2}+ \cdots+ a_{2,n}y^m
\end{equation*}
\noindent is a generic form of degree $m$ with $n$ terms and 
$I=\langle f_1,f_2 \rangle = \langle f_1,r \rangle$. Thus, without loss of
generality for Gr\"obner basis calculations,
one can start with the following forms:
\begin{eqnarray}\label{f2}
f_1 &=& a_{1,1}x^n+ a_{1,2}x^{n-1}y+ \cdots+ a_{1,n}xy^{n-1}+ a_{1,n+1}y^{n}, \nonumber \\
f_2 &=& a_{2,1}x^{n-1}y^{\mu+1}+ a_{2,2}x^{n-2}y^{\mu+2}+ \cdots+ a_{2,n}y^m.
\end{eqnarray}

Construct the set $G =\{f_1,f_2,\dots,f_{n+1}\}$ from $f_1$ and $f_2$ by taking
$f_{t+2} = \overline{S(f_t,f_{t+1})}^{f_{t+1}}$ for $1 \leq t \leq n-1$, where
$S(f_t,f_{t+1})$ is the $S$-polynomial of $f_t$ and $f_{t+1}$ (see, e.g.,
\cite[\S 2.6, Definition 4]{CLO}).
Then it is not hard to see that
\begin{eqnarray}\label{G}
f_1 &=& a_{1,1}x^n+ a_{1,2}x^{n-1}y+ \cdots+ a_{1,n}xy^{n-1}+ a_{1,n+1}y^{n}, \nonumber \\
f_2 &=& a_{2,1}x^{n-1}y^{\mu+1}+ a_{2,2}x^{n-2}y^{\mu+2}+ \cdots+  a_{2,n}y^m, \nonumber \\
f_3 &=& a_{3,1}x^{n-2}y^{\mu+3}+ a_{3,2}x^{n-3}y^{\mu+4}+ \cdots+  a_{3,n-1}y^{m+1}, \nonumber \\
    &\vdots&					\nonumber \\
f_t &=& a_{t,1}x^{n-(t-1)}y^{\mu+(2t-3)}+ \cdots+  a_{t,n-(t-2)}y^{m+(t-2)}, \\
 &\vdots&  \nonumber \\ 
f_n &=& a_{n,1}xy^{\mu+(2n-3)}+ \cdots+  a_{n,2}y^{m+(n-2)}, \nonumber \\
f_{n+1} &=&  a_{n+1,1}y^{\mu+(2n-1)}, \nonumber 
\end{eqnarray}

where, for $1 \leq  t \leq n-1$,

\begin{equation*}
a_{t+2,i} = \left\{ \begin{array}{cc}
	(\frac{a_{t,i+2}}{a_{t,1}}-\frac{a_{t+1,i+2}}{a_{t+1,1}})- \frac{a_{t+1,i+1}}{a_{t+1,1}}
		(\frac{a_{t,2}}{a_{t,1}}-\frac{a_{t+1,2}}{a_{t+1,1}}) & 1 \leq i \leq n-t-1, \\
	\frac{a_{t,n-t+2}}{a_{t,1}}-\frac{a_{t+1,n-t+1}}{a_{t+1,1}}
		(\frac{a_{t,2}}{a_{t,1}}-\frac{a_{t+1,2}}{a_{t+1,1}}) & i = n-t.
	\end{array}\right.
\end{equation*}
 
All the coefficients are nonzero, since the forms we started with are generic.
The following theorem shows that $G$ is a Gr\"obner basis for the ideal $I$.

\smallskip
\begin{theorem}\label{thm1}
Let $I= \langle f_1,f_2 \rangle$, where $f_1$ and $f_2$ as in (\ref{f2}), and let $G$
be the set in (\ref{G}). Then G is a Gr\"obner basis for the ideal $I$ with respect
to the graded reverse lexicographic order, with $x > y$.
\end{theorem}

\begin{proof}
By using \cite[\S 2.9 Theorem 9]{CLO} and the fact that 
$f_{t+2} = \overline{S(f_t,f_{t+1})}^{f_{t+1}}$ for $1 \leq t \leq n-1$, 
it is enough to
show that the set of syzygies 
$\mathcal{S} = \{S_{t,t+1} \: : \: 1 \leq t \leq n-1 \}$ forms
a homogeneous basis for the set of all syzygies 
$\{S_{i,j} \: : \: 1 \leq i < j \leq n+1 \}$ among the elements of $G$.
 
To simplify the calculations,
we are going to make all the polynomials of $G$  monic by dividing
each one by its leading coefficient; let $b_{s,t} = \frac{a_{s,t}}{a_{s,1}}$.
So without loss of generality, we will assume that all $f_i$ are monic.
It is easy to see that $S_{i,i+1}$ is a homogeneous syzygy. In fact,
\begin{equation}
S_{i,i+1} = \left\{ \begin{array}{cc}
		(y^{\mu+1},-x,0,\dots,0) & i = 1, \\
	(0,\dots,\underbrace{y^2}_{i-\text{th pos.}},-x,0,\dots,0) &  1 < i \leq n.
		\end{array} \right.
\end{equation}
Claim: For all $1 \leq i \leq n$ and $1 \leq t \leq n-i+1$, 
$S_{i,i+t}$ is generated by elements of $\mathcal{S}$. In fact
\begin{eqnarray}\label{6}
S_{i,i+t} &=& y^{2t-2}S_{i,i+1}+xy^{2t-4}S_{i+1,i+2}+\cdots+x^{t-1}S_{i+t-1,i+t} \nonumber\\
	  &=& \sum_{j=0}^{t-1} x^jy^{2(t-1-j)}S_{i+j,i+j+1}. 
\end{eqnarray}
Let $1 \leq i \leq n$.  The proof will proceed by induction on $t$. For $t=1$, 
$S_{i,i+t} = S_{i,i+1} \in \mathcal{S}$. Now assume the
claim is true for $t$, i.e.,
$S_{i,i+t}$ is a combination of elements of $\mathcal{S}$.

We need to show (\ref{6}) for $t+1$. 
Hence, we need to show that $S_{i,i+t+1}$ can
be written as a combination of elements of $\mathcal{S}$. For $i=1$,
\begin{eqnarray}
S_{1,1+(t+1)} &=& (y^{\mu+2(t+1)-1},0,\dots,0,\underbrace{-x^{t+1}}_{t+2-\text{th pos.}},
	0,\dots,0) \nonumber \\
	&=& y^2(y^{\mu+2t-1},0,\dots,0,\underbrace{-x^t}_{t+1-\text{th pos.}},0,\dots,0)+ \nonumber \\
	& &	x^t(0,\dots,0,y^2,\underbrace{-x}_{t+2-\text{th pos.}},0,\dots,0) \nonumber \\
	&=& y^2S_{1,1+t}+x^tS_{1+t,(1+t)+1}. \nonumber
\end{eqnarray}
And for $ i > 1$,

\begin{eqnarray}
S_{i,i+(t+1)} &=& (0,\dots,0,\underbrace{y^{2(t+1)}}_{i-\text{th pos.}},
			0,\dots,0\underbrace{-x^{t+1}}_{i+t+1-\text{th pos.}},0,\dots,0) \nonumber \\
	&=& y^2(0,\dots,0,\underbrace{y^{2t}}_{i-\text{th pos.}},0,\dots,0,
		\underbrace{-x^{t}}_{i+t-\text{th pos.}},0,\dots,0) +  \nonumber \\
	& & x^t(0,\dots,0,\underbrace{y^2}_{i+t-\text{th pos.}},-x,0,\dots,0) \nonumber \\
	&=& y^2S_{i,i+t+1}+x^tS_{i+t,(i+t)+1}. \nonumber
\end{eqnarray}
Thus, the set  $\mathcal{S}$ is a homogeneous basis of the set of all syzygies on $G$. 
Therefore, by \cite[\S 2.9 Theorem 9]{CLO} and the way that we constructed $G$,
$G$ is a Gr\"obner basis for the ideal $I$.
\end{proof}

\smallskip 
\begin{cor}\label{cor1}
Let $I$, $G$ and $R$ be 
as above, then the initial ideal of $I$ with respect to 
the graded reverse lexicographic order is: 
\begin{equation}\label{J}
J = \langle x^n,x^{n-1}y^{\mu+1},x^{n-2}y^{\mu+3}, \dots, xy^{\mu+2n-3},y^{\mu+2n-1}\rangle.
\end{equation}
\end{cor}

\begin{proof}
The proof follows from the fact that $G$ is a  
Gr\"obner basis for the ideal $I$ 
with respect to the graded reverse lexicographic order.
Hence the initial ideal of $I$ with respect to this order is
\begin{eqnarray*}
in(I) &=& \langle LM(f) : f \in G \rangle \\ 
     &=&  \langle x^n,x^{n-1}y^{\mu+1},x^{n-2}y^{\mu+3}, \dots, xy^{\mu+2n-3},y^{\mu+2n-1}\rangle.
\end{eqnarray*} 
\end{proof}

\smallskip 
\begin{lemma}\label{lem1}
The monomial ideal $J$ in (\ref{J}) is a reverse lexicographic ideal.
\end{lemma}   

\begin{proof}
Let $x^\al y^\be \in J$ and let  $x^sy^t$ such that 
$s+t = \al+\be$ and $s \geq \al$.  We
need to show that $x^sy^t \in J$. Since  $x^\al y^\be \in J$, there exist  
$\al_1, \be_1 \geq 0$ and an $i$ such that  
$x^\al y^\be = x^{\al_1}y^{\be_1}(x^{n-i}y^{\mu+2i-1})$. Thus $\al = \al_1+n-i$. Since
$s \geq \al$, there exists an $h \leq i$ such that
$s=\al_1+n-h$. But $s+t =\al+\be=\al_1+n+i+\be_1-\mu-1$, so 
\begin{eqnarray*}
t &=& \mu + h+i+\be_1-1 \\
	&\geq&  \mu + 2i+\be_1-1 \\
	&\geq& \mu +2i-1.
\end{eqnarray*} 
Therefore, $x^sy^t$ is divisible by $x^{n-i}y^{\mu+2i-1} \in J$ 
and hence $x^sy^t \in J$.  
\end{proof}

\smallskip 
\begin{figure} 
\begin{center} 
{ 
\leavevmode 
\hfil 
\epsfysize=3.5truein\epsfbox{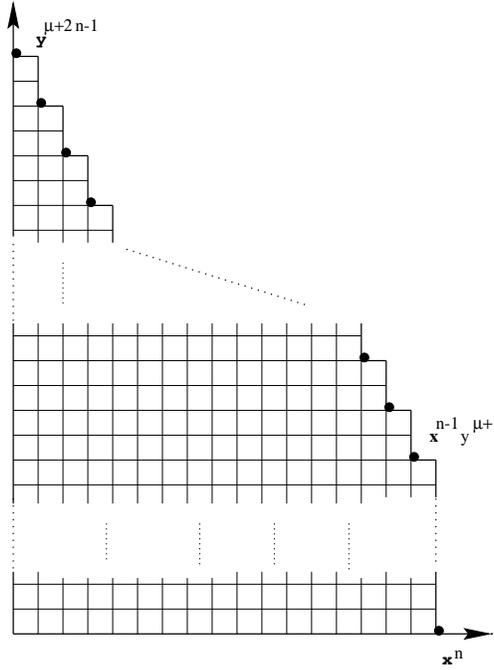} 
\hfil
} 
\caption{The graph  of the staircase of the ideal $J$ in Corollary \ref{cor1}. }
\label{fig1} 
\end{center}
\end{figure}

\smallskip 
\noindent Now we combine Corollary \ref{thm1} and Lemma \ref{lem1} to get
 a proof of the main theorem.

\smallskip 
\begin{theorem}[Moreno-Soc{\'\i}as Conjecture for two variables]
Let $R=K[x,y]$ be the polynomial ring  over an infinite field $K$. 
Let $d_1, d_2 \in \nn$
and let $f_1 \in R_{d_1}, f_2 \in  R_{d_2}$ be generic forms 
generating a generic ideal $I$. 
Let $J= in(I)$ be the initial ideal of $I$ 
with respect to the graded reverse lexicographic order.
Then $J$ is weakly reverse lexicographic.
\end{theorem}

\smallskip
\noindent As mentioned before, the main idea for the proof was to recognize 
the formula for the initial ideal from computations. We were not
able to do enough computations for higher $n$. 
Moreover, even if one could carry out enough calculations to guess 
a formula like (\ref{J}), we believe that it
would be hard to prove it: one needs to get a Gr\"obner 
basis like (\ref{G}), which will be hard to find 
using the same methods as for $n=2$.


\end{document}